\documentclass[leqno,amsmath,twocolumn]{article}
\usepackage{xspace}
\usepackage[english,italian]{babel}
\usepackage[latin1]{inputenc}
\usepackage[T1]{fontenc}
\usepackage{hyperref}
\usepackage{graphicx}
\usepackage{xcolor}
\def\site{https://wims.matapp.unimib.it/wims/wims.cgi}
\title{WIMS, a community of teachers, developers and users}
\hypersetup{colorlinks=true,urlcolor=blue,breaklinks=true}

\usepackage{subcaption}

\usepackage[noblocks]{authblk}

\author[*]{Marina Cazzola}
\author[**]{Sophie Lemaire}
\author[**]{Bernadette Perrin-Riou}
\affil[*]{Dipartimento di Matematica e Applicazioni, Università degli Studi di Milano-Bicocca, Milano, Italia.}
\affil[**]{Université Paris-Saclay, CNRS, Laboratoire de mathématiques d'Orsay, 91405, Orsay, France.}
\makeatletter
\@ifclasswith{article}{twocolumn}{%
\usepackage{geometry}
\geometry{
  paperwidth=8.25in,
  paperheight=10.75in,
  margin=1.4cm,
  includefoot
}
\setlength{\columnsep}{0.6cm}
}
\makeatother

\begin{document}

\selectlanguage{english}

\maketitle

\section{Introduction}

The so called ``new technologies'' are nowadays pervasive in every aspect of our lives, and their use in communication
of science and education is an established need.
The aims for such an use can be various. On the one hand, in the spirit of improving teaching, computers can be used to improve the standards of the materials that are offered to the students, as much as it is allowed by the new media. On a completely different point of view, the use of computerized tests can have the effect of reducing teachers' load of work, so to save energies that then can be directed towards tasks more significant than just ``marking''.
Furthermore the possibility to build personalized environments could put the students in the position to work at their own pace and develop self directed learning skills.

Having been engaged in the development of WIMS (WWW Interactive Multipurpose Server) since almost two decades, we wish to contribute with an analysis of our experiences about the benefits that these technologies can bring along.

\section{About WIMS}
The WIMS system started as an individual project by XIAO Gang, based at Département de mathématiques at Université de Nice - Sophia Antipolis (France), and had been made available to the public domain in 1998~\cite{xiao01:wims}. Xiao's project aimed at ``a systematic approach for providing Internet-accessible mathematical computations'' (\cite{xiao01:wims}).
It included, since the beginning, basic LMS-like functionalities: ``Virtual classes'' for teacher to dispense learning materials to the students and monitor their achievements, a grading system, a forum, \dots{} (for a complete description see~\cite{guerimand04:wims}).
Nowadays, under the name WIMS you can find a network of servers sharing interactive resources at many levels in various subjects (not only mathematics but also biology, chemistry, economics, languages, physics \dots{}), most of which have automatic correction and marking of user input.
The LMS functionalities in WIMS significantly evolved since the beginning (better user management, better feedback for users' inputs, e.g. through opportune links to specific ``Documents'' compiled by the teachers).
WIMS also provides a secure exam mode.

The characteristic for which today WIMS still stands out from the other currently available LMS is the built in possibility of extensively using random parameters in the coding of the learning objects\footnote{Random parameters can also be used in the editing of documents.} and the capability of interacting with softwares of the most different kinds.

These two features combined allow for the creation of complex and engaging exercises. A well designed use of random parameters enable WIMS to provide a ``virtually infinite'' set of copies of each single activity.

\subsection{A first approach to WIMS}

\begin{figure}[htbp]
  \centering
  \begin{subfigure}[h]{\linewidth}
  \centering
  \fbox{\includegraphics[width=0.8\linewidth]{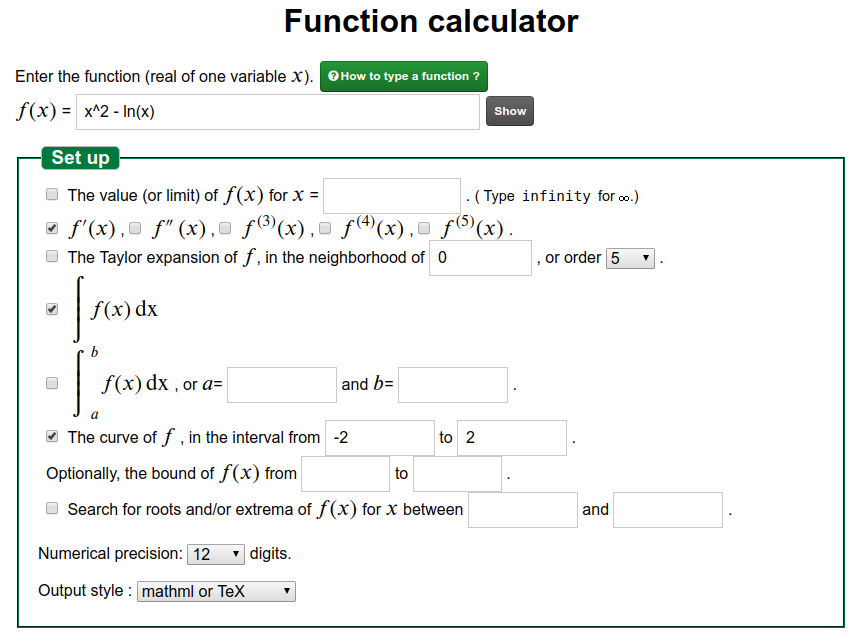}}
  \caption{tool/analysis/function.en}
  \label{fig:function.en}
  \end{subfigure}
  \begin{subfigure}[h]{\linewidth}
  \centering
  \fbox{\includegraphics[width=0.8\linewidth]{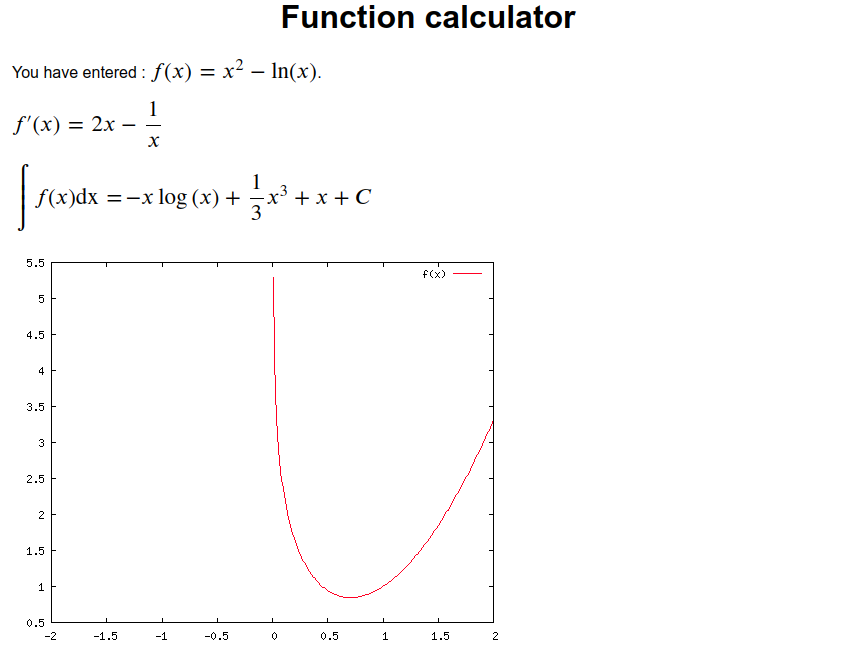}}
  \caption{WIMS' answer}
  \label{fig:function.en:rep}
  \end{subfigure}
  \caption{WIMS tool}
\end{figure}

The simplest use of WIMS in mathematics is to let it act as a user-friendly web interface for asking for computations. For example, the tool ``Function calculator''~\cite{tool:function} interacts with public domain softwares as Maxima\footnote{\url{http://maxima.sourceforge.net/}}, PARI/GP\footnote{\url{http://pari.math.u-bordeaux.fr/}} and Gnuplot\footnote{\url{http://www.gnuplot.info/}}. The user inserts the query in the mask shown in Figure~\ref{fig:function.en} and gets the answer shown in Figure~\ref{fig:function.en:rep}: it is up to WIMS to choose, among all the possibilities, the software more suitable for each task.\footnote{%
The tool is freely available on any WIMS server, e.g. you can access
and test it at
\url{\site?module=tool/analysis/function.en}.}

\begin{figure}[htbp]
  \centering
  \begin{subfigure}[h]{\linewidth}
  \centering
  \fbox{\includegraphics[width=0.9\linewidth]{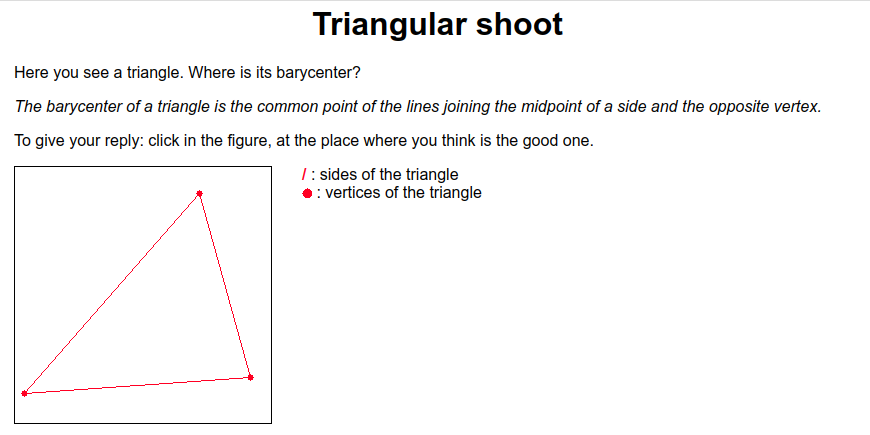}}
  \caption{Barycenter I}
  \label{fig:trishoot:1}
  \end{subfigure}
  \begin{subfigure}[h]{\linewidth}
  \centering
  \fbox{\includegraphics[width=0.9\linewidth]{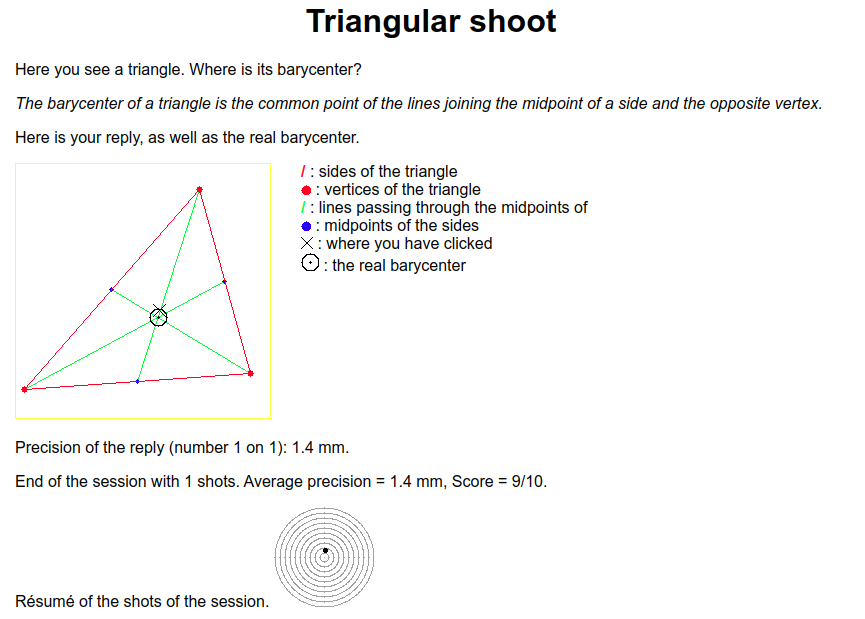}}
  \caption{Barycenter: analysis of the answer}
  \label{fig:trishoot:2}
  \end{subfigure}
  \caption{Barycenter}
\end{figure}

\begin{figure}
  \centering
  \fbox{\includegraphics[trim=0 0 0 70,clip,width=0.9\linewidth]{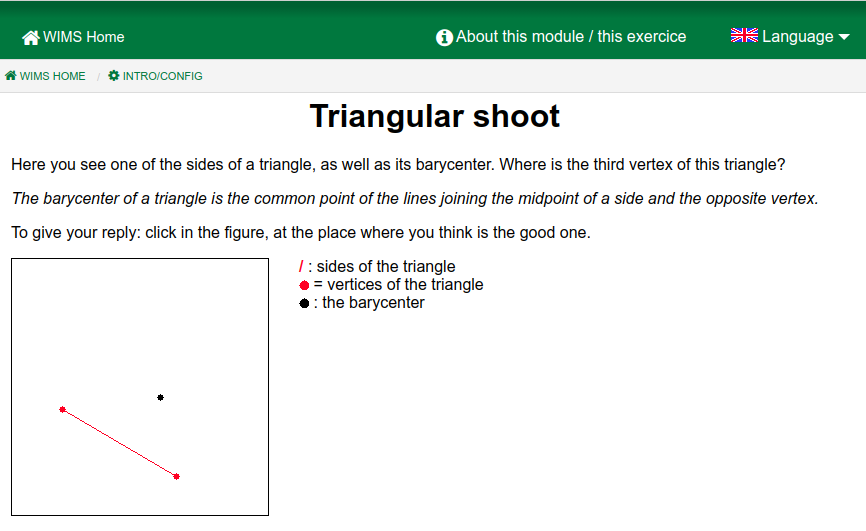}}
  \caption{Barycenter: inverse problem}
  \label{fig:trishoot:3}
\end{figure}

If we focus on the more complex task of managing interactive activities and effectively analyzing user answers,
one of the strengths of WIMS is its ability to interact with representations and drawings. For example, the activity ``Triangular shoot''~\cite{trishoot} allows the users to familiarise with the centers of a triangle (barycenter, orthocenter, incenter and circumcenter). Every time the exercise is selected, a different randomly generated triangle is proposed. In the easier version of the exercise, a triangle is given and the user has to click on one of the centers (see Figure~\ref{fig:trishoot:1}). The user's reply is evaluated and marked (see Figure~\ref{fig:trishoot:2}).
To better experience with such concepts, ``Triangular shoot'' also includes the reverse problem (see Figure~\ref{fig:trishoot:3}): given two vertices of the triangle and its barycenter, can you identify the missing vertex of the triangle?

\subsection{The structure of WIMS}
WIMS has a modular structure (``module'' is the term used by WIMS for any of its ``units''). A standard installation of the latest version of WIMS includes about 50 administrative modules and about 1800 learning modules.
Whenever a new task is required, a new administrative module can be added to the system to provide such functionality.
The learning modules include what WIMS calls ``tools'' (e.g. the already described ``Function calculator'') as well as learning objects contributed by the teachers themselves (exercises, virtual classes, manuals \dots) designed according to their individual needs and then shared with the whole community.

As mathematics is concerned, WIMS is fully \LaTeX{} aware, and it can interact with the already mentioned Maxima, PARI/GP and Gnuplot. It can also interact with GAP\footnote{\url{http://www.gap-system.org}} (e.g. see ``OEF permutation''~\cite{oefperm}), Povray\footnote{\url{http://www.povray.org/}} (e.g. see ``Polyray''~\cite{polyray}), Graphviz\footnote{\url{http://www.graphviz.org/}} (e.g. see ``Graph drawing tool''~\cite{graphviz}), GeoGebra\footnote{\url{https://www.geogebra.org/}} (e.g. see ``OEF GeoGebra''~\cite{oefgeogebra}), Octave\footnote{\url{https://www.gnu.org/software/octave/}} (e.g. see ``Statistical tables''~\cite{statisticaltables}).
More recently the core of WIMS has been updated in order to interact with javascript applets.
The drawing capabilities of JSXGraph\footnote{\url{http://jsxgraph.org}}
and of HTML5 Canvas are now integrated into WIMS allowing much more
interactivity from the users.
Other applets have been interfaced allowing for the creation of exercises in other subjects. For example, JSME\footnote{\url{https://peter-ertl.com/jsme}} and JSmol\footnote{\url{http://jmol.sourceforge.net}} are used in exercises on molecular representations in the module ``Amino acids''~\cite{oefaminoac}. JSmol is also used to draw  polyhedra (see the online tool ``Convex Polyhedra''~\cite{polyhedra}).

\begin{figure}[htbp]
  \centering
  \fbox{\includegraphics[trim=0 110 0 800,clip,width=0.9\linewidth]{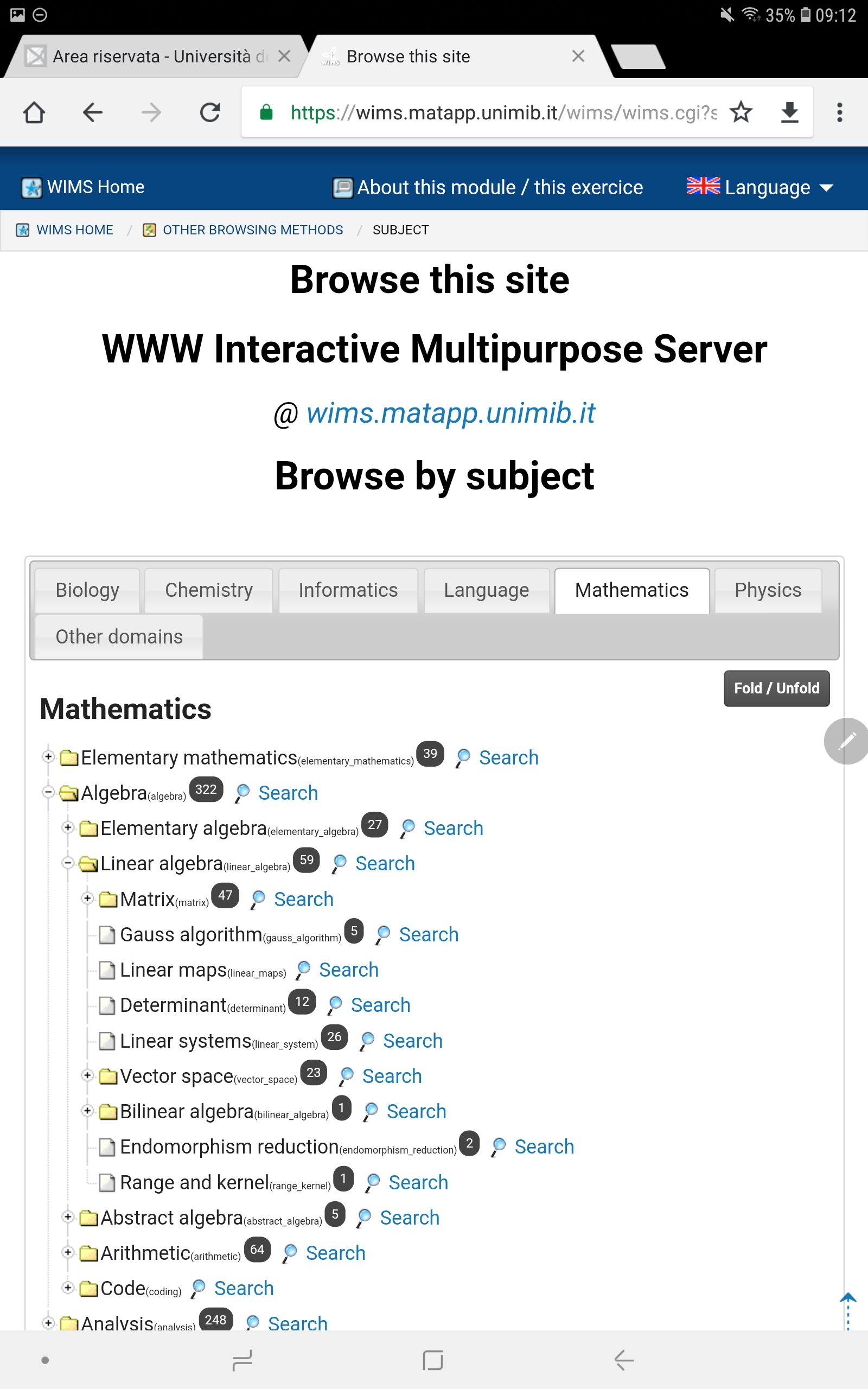}}
  \caption{Browse by subject}
  \label{fig:subject}
\end{figure}

All of WIMS' learning modules are classified by the authors via keywords according to the topics covered, and a search engine is provided. The ever expanding tree of keywords can be navigated via the ``Browse by subject'' module (see Figure~\ref{fig:subject}).
Classifications according to school levels (from pre-school to university postdoc) are also provided.

\section{The community}

Started as a one-man project, WIMS owes its growth to the large community that Xiao was able to gather around his work.
The project was published in the public domain, giving open access to developers and to translators (the core of WIMS is available in French, English, Italian, Dutch, Chinese, Slovene, Spanish and Catalan).

Xiao also called for the cooperation of teachers, promoting for them training seminars on the programming of exercises.
From its very beginning, the development of WIMS relies on contacts between the developers and the community of teachers, who constantly report bugs and give ideas for improvements.
Since 2006 periodical ``Colloquia'' are held (e.g. see~\cite{coll-2-nice:wims}), in order to discuss progresses in the use of WIMS and possible development strategies.

In July 2007, WIMSEDU (\emph{Enseignants, développeurs et utilisateurs de WIMS})\footnote{\url{https://wimsedu.info}} was founded, an association with the aim of strengthening this community: share skills, act as interlocutors for WIMS' users towards users of other software or institutions,
discuss and influence the directions of its development by supporting the chosen projects,
create a dynamic collaborative work through its website, identify,
support and inform about training actions.
WIMSEDU carries on the tradition of the ``Colloques'' (Nice 2006, Nice 2007, Bordeaux 2010, Rennes 2012, Dunkerque 2014, Annecy 2016, Paris Orsay 2018, Amiens 2020) and the meeting of different experiences continues (mathematicians and non-mathematicians, developers and users, researchers and teachers from any educational level). In the Colloquia experiences of use are confronted, technical difficulties are discussed and innovations are proposed.
The association pursues its goals of acting as the glue of this community also by directly managing a server of distribution of learning modules through which the authors can share the exercises and documents they have designed for their students so that they are published on all of WIMS servers.
The actual development of the core of WIMS is carried out by a volunteers development team, whose activities are intertwined with WIMSEDU actions. Many of the team members are university staff, but no academic institution is directly involved, apart from hosting main servers in the WIMS network (servers generally used for institutional didactics). An important share of the work is still carried out by teachers (of any school levels: from primary school teachers to university professors).
The source code of the core of the system, in a ``public domain'' spirit, is hosted by the French Public Interest Group RENATER\footnote{French National Telecommunication Network for Technology, Education and Research}, which in this way provides an essential suport for WIMS~\footnote{The project page can be reached at the following url \url{https://sourcesup.renater.fr/projects/wimsdev/}}.

\section{WIMS for mathematics education}

WIMS can be used for mathematics teaching as it allows to provide students with a big amount of engaging exercises. There are many descriptions of experiences with WIMS for teaching (e.g. see~\cite{vanderb06}, \cite{cazzola11:wimstd54}, \cite{ducrocq2010utilisation}, \cite{kuzman2007interactive}, \cite{ramage04:la}, \cite{reyssat13:enseigwims}, \cite{kobylanski2019wims}, \cite{cazzola2019:wims}), showing (at different school levels, in different class situations, in different subjects) that WIMS-based exercises, with automatic marking, can motivate students and keep alive their attention.
Moreover, WIMS can be used for exams, throught special functions that guarantee an appropriate level of security and control.

\subsection{Mathematical problem solving}

\begin{figure}[htbp]
  \centering
  \begin{subfigure}[h]{\linewidth}
     \centering
     \fbox{\includegraphics[width=0.97\linewidth]{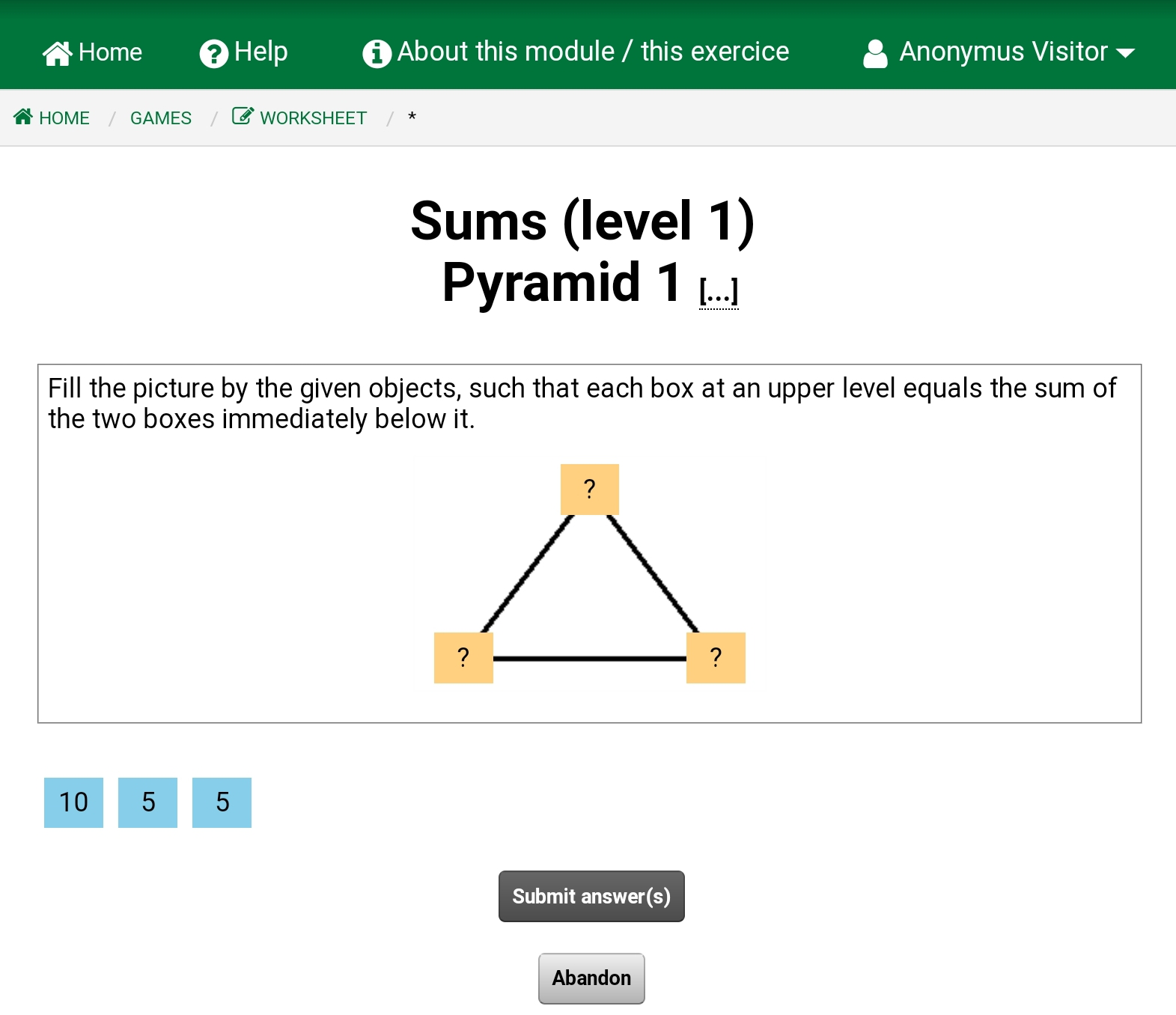}}
     \caption{Number pyramids I}
     \label{fig:sumsbase}
  \end{subfigure}
  \hspace{5pt}
  \begin{subfigure}[h]{\linewidth}
     \centering
     \fbox{\includegraphics[width=0.97\linewidth]{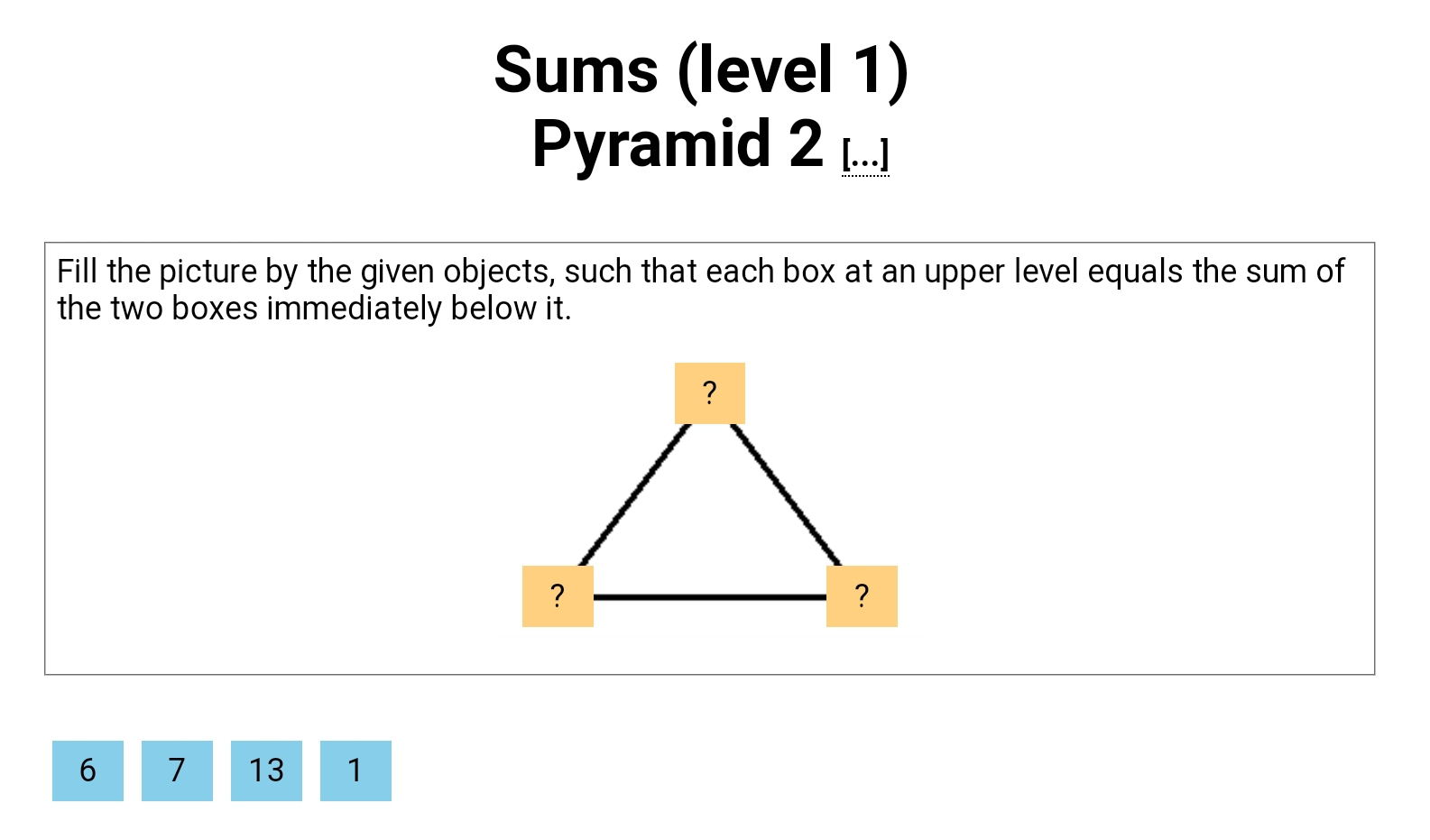}}
     \caption{Number pyramids II}
     \label{fig:forcesums}
  \end{subfigure}
  \caption{Number pyramids}
\end{figure}

We believe that problem solving \emph{à la} Polya is a crucial step for a real understanding of mathematics. To gain a ``taste'' for mathematics, pupils should be lead to actively work on a difficult problem and to develop a viable solution for such a problem, given the necessary amount of time~\cite{mathdiscovery}.
The use of ``automatic'' exercises might conflict with this idea. By beeing given a large amount of similar exercises, students might get the idea that the goal of mathematics is to find the \emph{fastest} strategy to in order to get a ``computer approved'' answer, strategy that not always corresponds to a real understanding of the subject. This attitude is even stronger if WIMS exercises are used for exams.
It is indeed a difficult equilibrium and a particular effort must be put on the design and selection of the learning objects and on the monitoring of students' activity on these exercises.
For example, WIMS proposes standard ``Number pyramids'' drills for practicing sums (e.g. at a primary school level) as the one shown in Figure~\ref{fig:sumsbase}.
If compared to a ``paper and pencil'' corresponding version of this exercise, the advantage of this computer based version is that the student have an immediate feedback and can repeat the exercise as many times as wanted.
It is easy to predict, and in fact we could observe the phenomenon, that this exercise can lead the students to develop ``improper'' strategies: after a few tries the students can understand that there is actually no need of calculating any sum, but a ``computer validated'' solution of the exercise can be obtained very quickly just putting the biggest number above and the two smaller numbers below.
This way of finding the answer for the exercise can be considered ``wrong'' if the aim of the teacher is having pupils work with computation of sums.
The teacher should monitor students' work and, if needed, should act so to avoid students to rely on this unwanted strategy, and lead them to proceed a step forward to tackle more advanced tasks. A simple way to gain this is to switch to a version of the same exercise which proposes more option for the numbers to be dragged into the pyramid, as the one shown in Figure~\ref{fig:forcesums}. The presence of an extra number to be set aside forces the user to actually ``do'' the sums (or to find a more creative strategy!).
The example in Figure~\ref{fig:forcesums} also allows us to show that in WIMS you can design exercises with multiple correct solutions, thus countering the formation of one of the most deleterious students' belief: ``Mathematics problems have one and only one right answer'' (cfr. \cite[p.~359]{schoenfeld92}).

\begin{figure}[htbp]
  \centering
     \fbox{\includegraphics[width=0.97\linewidth]{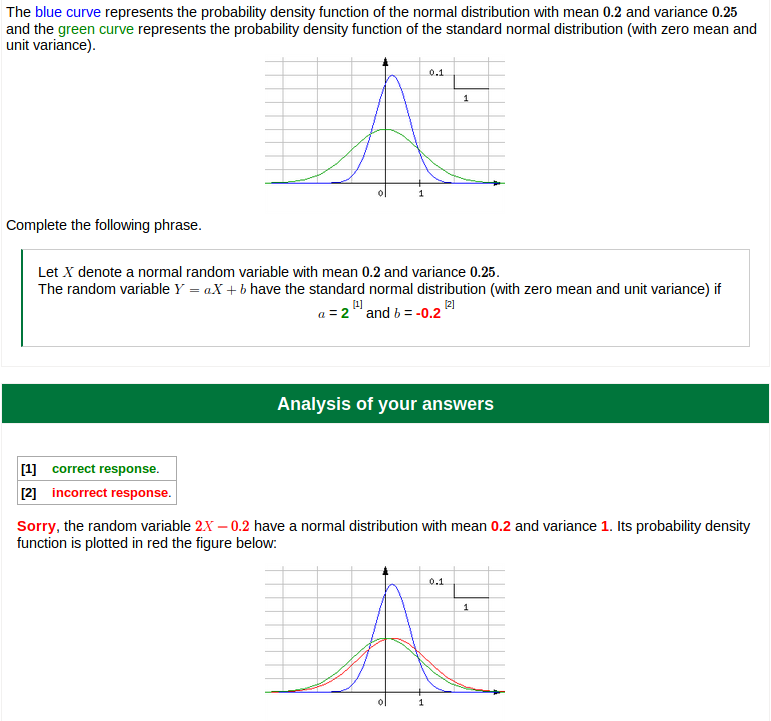}}
     \caption{Statement and feedback in exercise ``A standard normal random variable''}
     \label{fig:normaldistr}
\end{figure}

Furthermore, in a problem solving perspective, presenting automatically a full solution for the exercise as soon as the student makes a mistake can be counterproductive, as this could prevent the actual understanding of the errors and interfere with the pursuit for an autonomous solution. WIMS allows for the creation of feedbacks that provide alerts for classical errors or hints to help students to persevere in their search. For example the module ``Probability distribution plots''~\cite{oefproba} contains an activity that shows how a feedback can indicate to the student an expected property that is not satisfied by his answer (see Figure~\ref{fig:normaldistr}).
Teachers that use WIMS with their students can interact with the exercises's author and suggest suitable feedbacks.

\subsection{Mathematicians' point of view}
Mathematicians can contribute to the development of WIMS by pointing out which are the aspects of the discipline that the learning objects should focus on. Going back to ``Triangular shoot'' we can notice that the author of the exercise (a mathematician!) did not want to test users on the question ``what is the barycenter of a triangle'', as the definition of such point is right in the statement shown to the user (as you can see in Figure~\ref{fig:trishoot:1}). Rather, the aim of the exercise is to use the definition and actively ``manipulate'' such notion.
Too often students are led to believe that most of mathematics consists of pure memorization of facts, and do not manage to get any real insight of the different facets of each concept.

\begin{figure}[htbp]
 \centering
 \fbox{\includegraphics[width=0.9\linewidth]{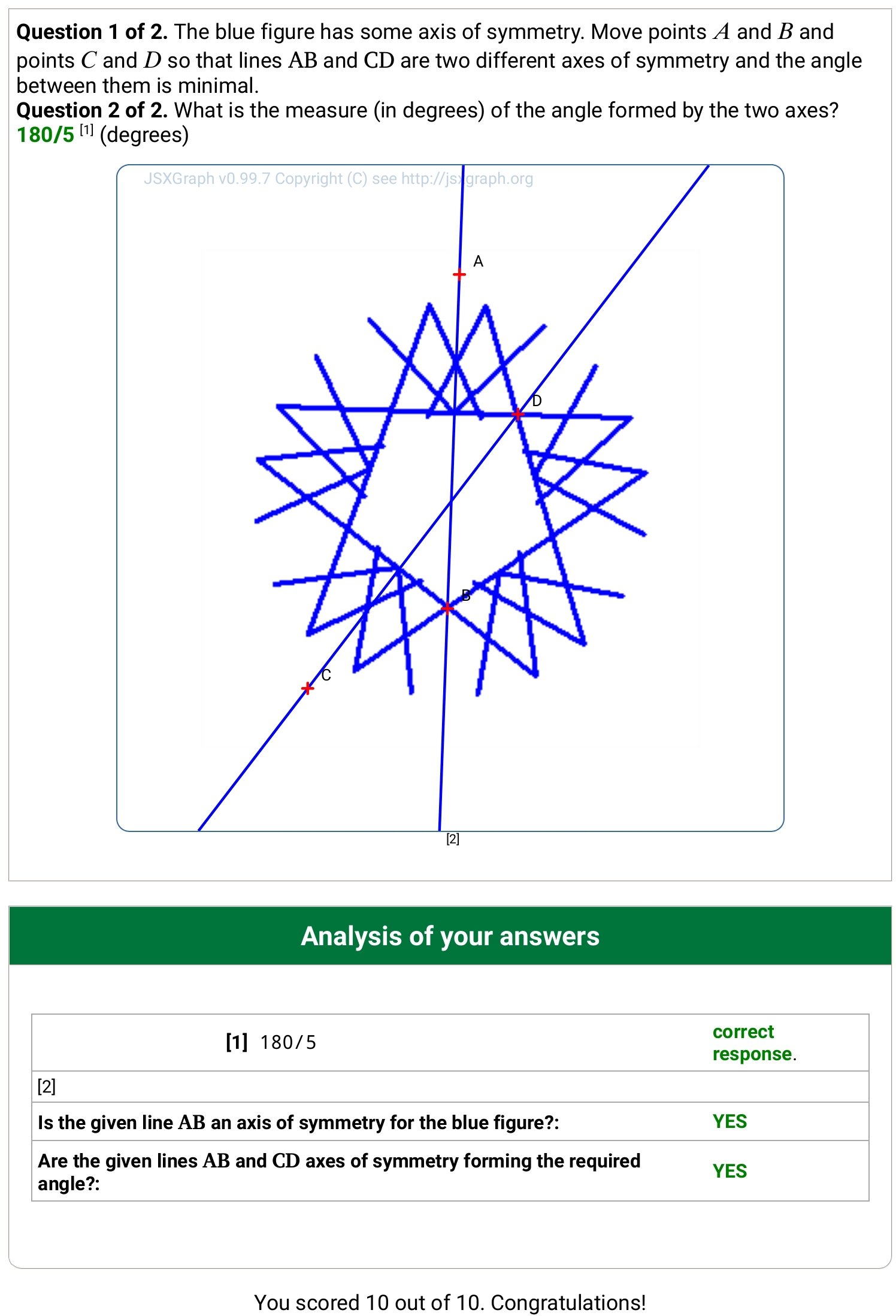}}
 \caption{Symmetry axes (level ``D'')}
 \label{fig:symm:5}
\end{figure}

A problem should be the mean to lead students grasp the theory underneath it.
An example in this direction is the series of activities on symmetry that are proposed: these exercises are built on different levels, to gradually lead learners to tackle more and more complex concepts.  E.g., in the module ``OEF Rosettes''~\cite{rosettes}, you can start (level ``A'') with the simple task of tracing any symmetry axis of a given figure, and later you can proceed to the more complex task (level ``D'', shown in Figure~\ref{fig:symm:5}) of finding the two symmetry axes that generate the whole \(*n\bullet\) kaleidoscopes rosette pattern (e.g. see~\cite{conway:magic}).

\begin{figure}[htpb]
  \centering
  \begin{subfigure}[b]{\linewidth}
  \centering
  \fbox{\includegraphics[width=0.9\linewidth]{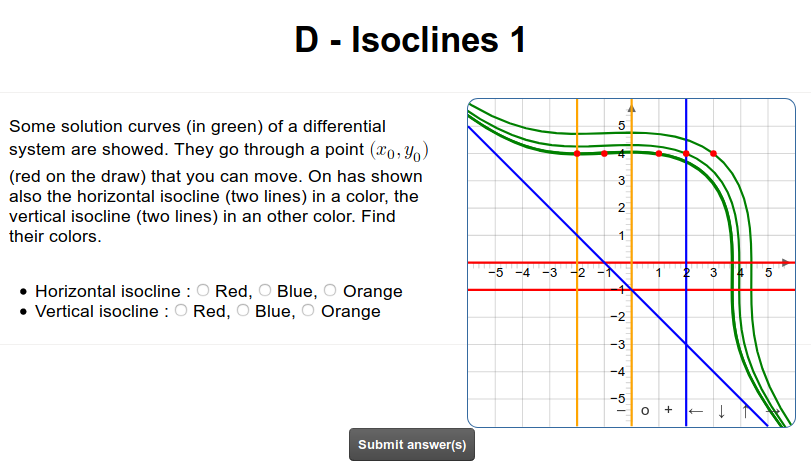}}
  \caption{Statement}
  \end{subfigure}
  \begin{subfigure}[b]{\linewidth}
  \centering
  \fbox{\parbox{0.9\linewidth}{
      \includegraphics[width=0.29\linewidth]{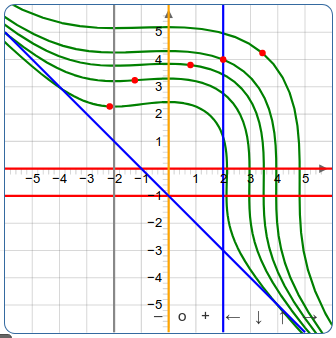}
      \includegraphics[width=0.29\linewidth]{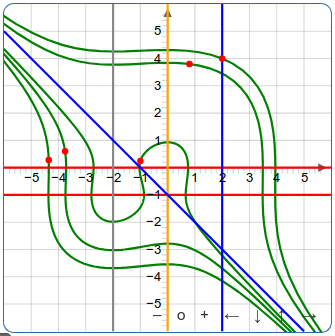}
      \includegraphics[width=0.29\linewidth]{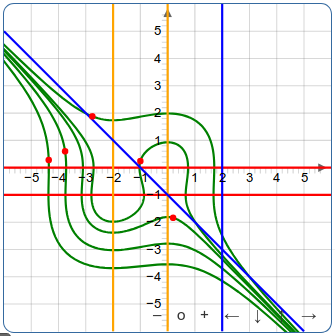}
    }}
  \caption{Manipulations by the user}
  \end{subfigure}
  \caption{Isoclines}
  \label{fig:isocl}
\end{figure}

\begin{figure}[htbp]
 \centering
 \fbox{\includegraphics[width=0.9\linewidth]{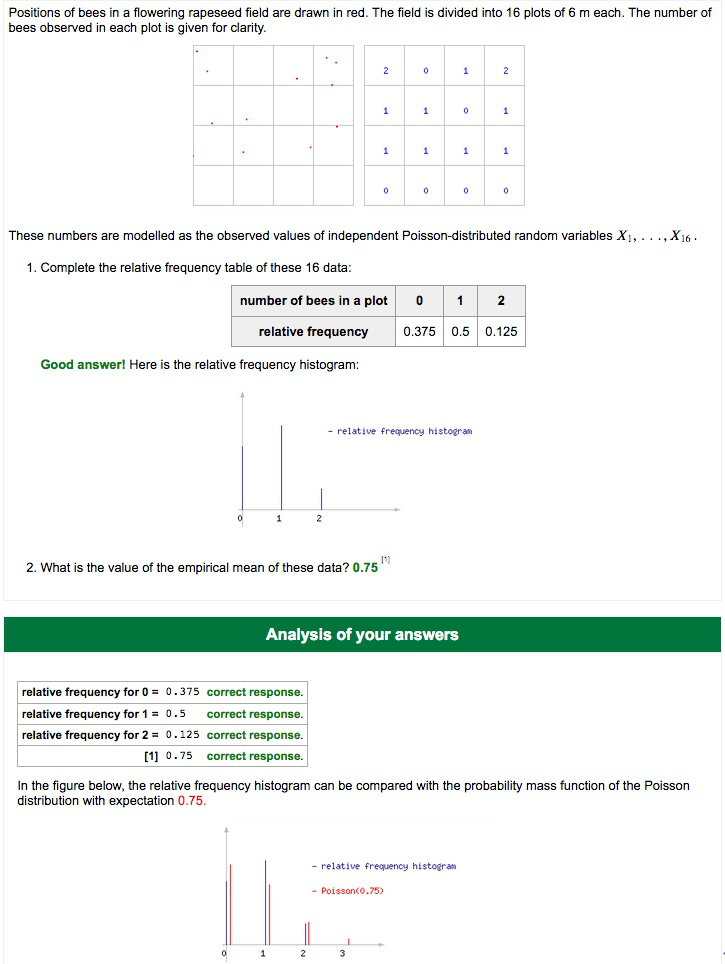}}
\caption{A statement of exercise ``Poisson distribution'' with a feedback when the student's answers are correct. }
 \label{fig:poisson}
\end{figure}

The idea of manipulating mathematical objects can be exploited in various fields of mathematics. The module ``{OEF} Graphics study of differential equations or simple differential systems''~\cite{grapheqdiff} contains an activity in which the user can interact with the solution curves of differential system, in order to gain information about its isoclines (see Figure~\ref{fig:isocl}).
A similar manipulation can be useful in other domains, such as probability, where many results are counter-intuitive for beginners.  Exercises in which simulated data and probability distributions are presented can help students to have a better understanding of what kind of events may appear frequently by chance and what a probability distribution represents. E.g. in the exercise ``Poisson distribution'' of the module ``Point estimation in statistics''~\cite{eststat}, the outcome of  a Poisson point process simulation on a square is drawn. The aim of the exercise is to find the empirical measure of the number of points in each cell of the given Cartesian grid and to estimate the intensity of this process (see Figure~\ref{fig:poisson}).
Once the student has given his answer, the relative frequencies of the number of points in each cell and the probability mass function of the Poisson distribution with the estimated parameter are shown in the same grouped bar chart through a feedback.
The fact that the two sets of bars can be quite different surprises some students.

Traditional problems can become a source of inspiration for creating learning objects.
For example, the problem of studying \emph{lattice polygons}, that is plane polygons with vertices with integer coordinates, has been widely explorered
and can give rise to nice results as Pick's theorem (e.g. see~\cite[p.~208-209]{coxeter1969}).
The module ``OEF Polygons on graph paper''~\cite{polyqq} offers a series of exercises
asking for the tracing of triangles and quadrilaterals with integer coordinates vertices.
The rationale behind the module is that exploring the variety of existing examples of
polygons opens the mind to a more advanced geometric vision, so the aim of the exercises
is to stimulate for the construction of non-stereotyped examples, by asking questions as:
\begin{itemize}
\item trace a quadrilateral with all equal sides and with no sides laid on the grid,
\item trace a quadrilateral with perpendicular and equal diagonals and with sides all of different lengths,
\item trace a quadrilateral with exactly two right angles and with no sides laid on the grid.
\end{itemize}
It might be interesting to point out that users' replies for such constructions are evaluated via basic 2D vector geometry notions, using PARI/GP: the coding of the exercise itself becomes a playground to see applied mathematics at work.

\begin{figure}[htbp]
  \centering
  \fbox{\includegraphics[width=0.9\linewidth]{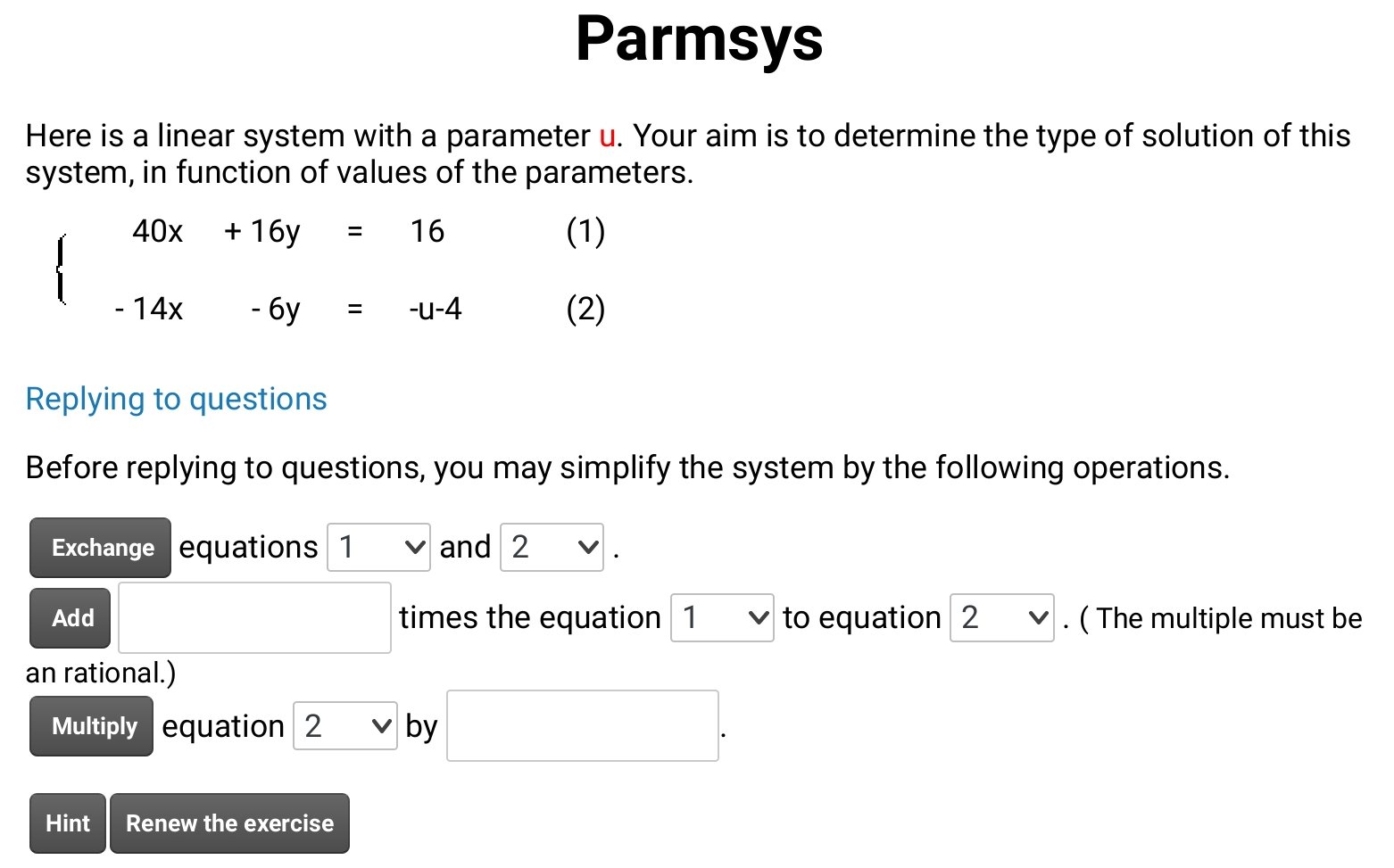}}
  \caption{Gauss' algorithm in linear algebra (step 1)}
  \label{fig:linalg:2}
\end{figure}

A computer based learning system can be used to introduce ``standard'' algorithms. In WIMS you can find, for instance, the module ``Parmsys''~\cite{parmsys} that leads the user through the steps of Gauss elimination method for reducing a linear system with parameters in order to determine whether it is solvable or not (Figure~\ref{fig:linalg:2}).
Activities of this kind, in which the task of actually making all the computations is left to the computer, allow for the students to focus on the problem itself: which are the most convenient steps to get to the triangular form for the system?
Once the algorithm has been assimilated, it is possible to propose modeling exercises for different applications. The module ``OEF linear systems''~\cite{oeflin} contains the exercise ``Alloy of 3 metals'' in which students are explicitly invited to use the tool ``Linear solver''~\cite{linsolver} so that they do not need to spend too much time on computations, but they can concentrate on the main difficulty of this type of exercises which is indeed to define the linear system they have to solve to answer the questions.

\begin{figure}[htbp]
  \centering
  \begin{subfigure}[b]{\linewidth}
  \centering
  \fbox{\includegraphics[width=0.9\linewidth]{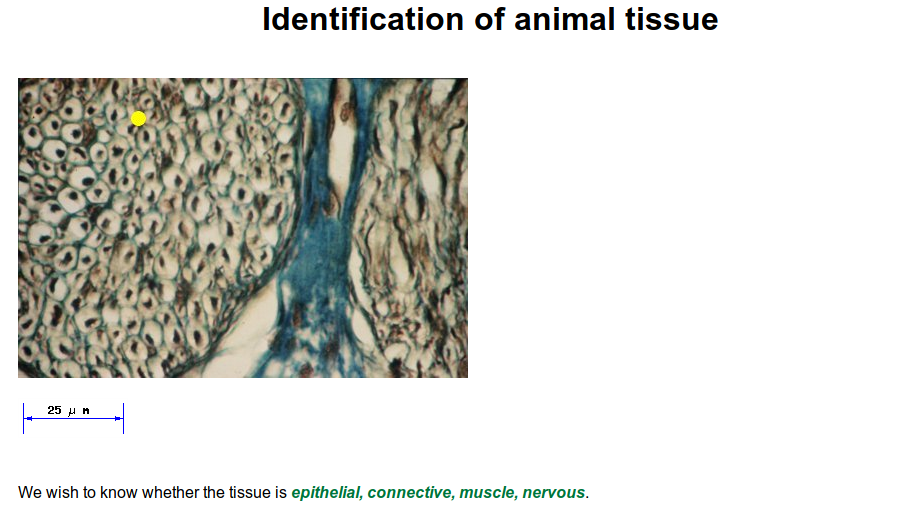}}
  \caption{Animal tissue I}
  \label{fig:bio1}
  \end{subfigure}
  \begin{subfigure}[b]{\linewidth}
  \centering
  \fbox{\includegraphics[width=0.9\linewidth]{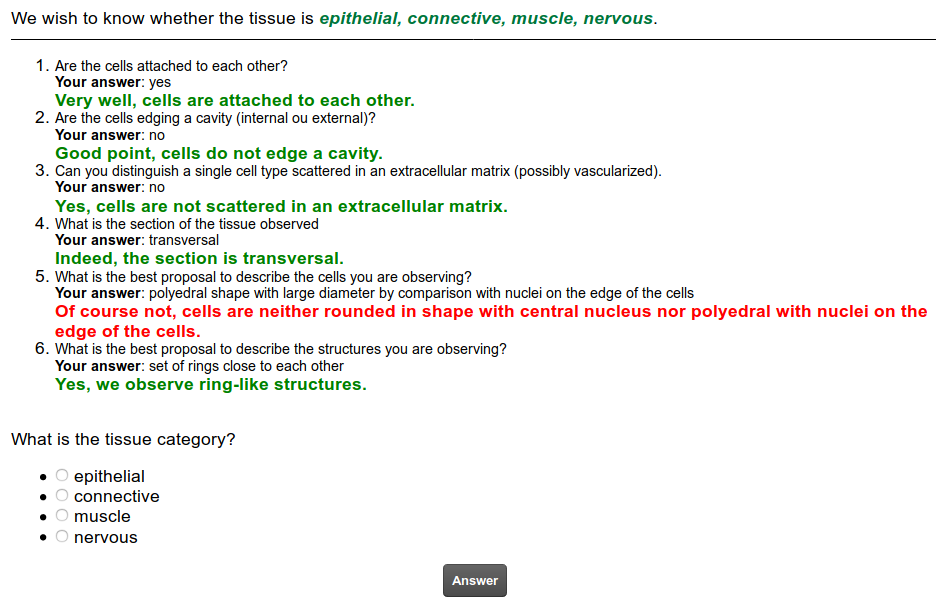}}
  \caption{Animal tissue II}
  \label{fig:bio2}
  \end{subfigure}
  \caption{Example in biology}
\end{figure}

Finally, and this time we exit the domain of mathematics and provide an example in biology, even more interactive activities can be build. In the module ``Identification of animal tissue''~\cite{dialbio}, a picture of a biological tissue to be identified is shown (see Figure~\ref{fig:bio1}), the user is guided through a standard analysis protocol with \emph{ad hoc} questions.
Users' answers are checked by the system one by one (see Figure~\ref{fig:bio2}), and, in case of errors, intermediate questions can be added so to allow for the student to rectify any mistake or to follow an alternave approach.
Only after any single reply has been given a feedback, the final question ``What is the tissue category?'' is asked. This example also allows us to see the idea of ``data modules'' at work: the set of images of tissues to be identified
is available to be shared with any other learning module on similar subjects.
\subsection{Teachers' point of view}

Through any WIMS server, teachers have access to a large variety of learning objects, on many different subjects, beyond mathematics. Available ``Open classes'' provides models for ready-to-use courses on various topics, so that teachers, with just one click, can create their own copy of the class and have access to all of the LMS functions available (selecting exercises, monitoring students' progresses, grading, \dots{}). Teachers have thus the option to use already existing materials ``as is'', but they also have the option to adapt anything to their own teaching style.
Moreover, they can create their own exercises. WIMS provides a somehow simplified programming language, called OEF (``Online Exercise Format''), that allows for the use of all the advanced features of WIMS and it is fully documented on line. WIMS also provides tools to create standard exercises, like multiple choice questions or ``true or false'', or ``fill in the blanks'' and so on.
Each WIMS module in the distribution is open source, so the teacher is encouraged to ``look inside'' any exercise and have an understanding on ``how it really works'', not merely accepting mysterious ``turnkey'' packages.

Teachers have the possibility to select among many different configurations for each single instance of any OEF exercise: if the author of the exercise provides a hint or a solution, it is up to the teacher to decide whether such hints or solution are actually shown to the students.
Teachers can write themselves alternative hints or solution to be shown to the students, using WIMS ``documents'' (in this case with no need to edit the exercise code).
Teachers can also setup the exercise so that in case of errors the student have to go through the very same exercise again: as WIMS exercises are highly randomized, usually you can expect a different version of the exercise any time you open it; with this option students have the chance to revise by themselves their resolution and empower their own self directed learning skills.

\begin{figure}[!h]
  \centering
  \fbox{\includegraphics[width=0.9\linewidth]{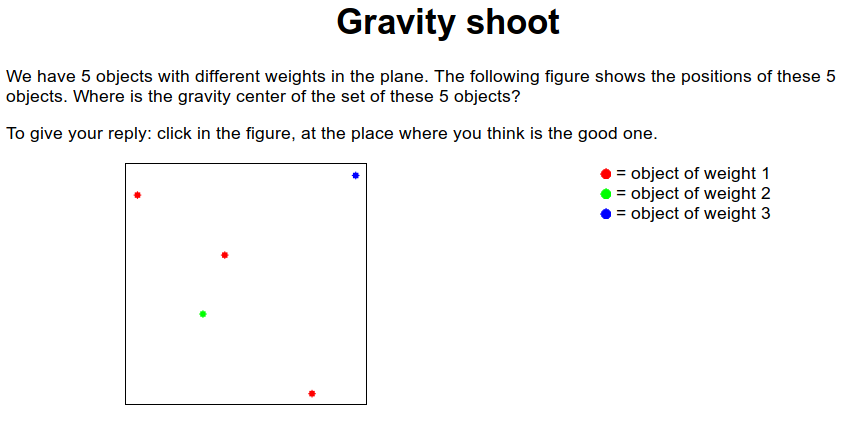}}
  \caption{Gravity}
  \label{fig:gravity:1}
\end{figure}

Finally, the fact that WIMS covers many subjects apart from mathematics allows teachers to have everything in one place and favors interdisciplinary links. For example, it allows for the connection of the idea of the Barycenter of a triangle seen in Figure~\ref{fig:trishoot:1} to the problem of finding the center of gravity of a system of weighed objects as in the module ``Gravity shoot'' \cite{gravshoot} shown in Figure~\ref{fig:gravity:1}.

\section{Conclusion}
We have tried to highlight the most significant aspects of the WIMS system and to bring our experience on the potential of using new technologies in teaching.
We hope to have encouraged the reader to visit the site, not letting themselves be restrained by a sometimes a bit intricate graphical aspect, but willing to really explore the numerous learning objects available.
We also hope to have encouraged participation in our community: there is a lot of work to be done in order to keep the WIMS system active and make the already available materials better usable.
A special effort in this direction concerns the translation into other languages of the huge resource of learning modules existing in French.
We also feel the need of building a better database of available resources and of improving examples of turnkey ``Open classes'', and, again, to translate modules and classes into many more languages.
Work has to be done to guarantee the operativity of modules in spite softwares obsolescence and in order to take advantage of technologies development: Xiao's modules, dating back as to 1998, are still fully functional.
And we constantly need to plan the future development of the system.

Once more we recall the questions raised in~\cite{cazzola13:ictmt}, that still remain open:
\begin{itemize}
\item Which are the most effective interactive activities, and can the
  strengths of WIMS help teachers offer such activities to their
  students?
\item Can WIMS be used to evaluate students' achievements?
\item Can WIMS be used to build personalized learning paths?
\item How to develop the technological knowledge of teachers in order
  to enable them to exploit the full potential of WIMS?
\item Can the capabilities of WIMS be embedded into other LMS, such
  as e.g. Moodle?
\end{itemize}
Any contributions from researchers in the subjects covered, from computer scientists and from researchers in education can support the development of WIMS in a direction that improves its effectiveness in communication of sciences and in teaching.

\bibliography{wimsmod,ho,wims}
\bibliographystyle{amsplainhyper}

\end{document}